\documentclass[12pt]{article}
 \usepackage[english]{babel}
 \usepackage{amssymb}
 \usepackage{amsmath}
 \usepackage{graphicx}
 \usepackage{bbm}
 \usepackage[right]{eurosym}
 \usepackage[labelfont={bf,sf},font={small},%
 labelsep=space]{caption}
 \usepackage{chngcntr}
 \usepackage{xcolor}
 
 \counterwithin{figure}{section}
 
 \def\vs{\vspace}

 \def\IR{\mathbb R}

 \def\exp{\mathrm{exp}}
 
 \def\ma{\mathcal}

\begin{document}

\begin{center}
{\Large\bf Restricted Log-Exp-Analytic Power Functions}

\vs{0.5cm}
Andre Opris
\end{center}

\vs{0.2cm}
{\small {\bf Abstract.} A preparation theorem for compositions of restricted log-exp-analytic functions and power functions of the form 
$$h: \mathbb{R} \to \mathbb{R}, x \mapsto \left\{\begin{array}{ll} x^r, & x > 0, \\
0, & \textnormal{ else, } \end{array}\right.$$
for $r \in \mathbb{R}$ is given. Consequently we obtain a parametric version of Tamm's theorem for this class of functions which is indeed a full generalisation of the parametric version of Tamm's theorem for $\mathbb{R}_{\textnormal{an}}^{\mathbb{R}}$-definable functions.

\vs{0.7cm}
{\Large\bf Introduction}

\vs{0.2cm}
In \cite{8} Opris gave the definition for restricted log-exp-analytic functions. These are $\mathbb{R}_{\textnormal{an,exp}}$-definable functions which are compositions of log-analytic functions and exponentials of functions which are locally bounded where $\mathbb{R}_{\textnormal{an,exp}}$ is the structure generated by all restricted analytic functions and the global exponential function (see \cite{1}). A log-analytic function is piecewise given by compositions from either side of globally subanalytic functions and the global logarithm (see \cite{5}, \cite{6} and \cite{9} for the formal definition and elementary properties of log-analytic functions).

\vs{0.4cm}
{\bf Example}

\vs{0.1cm}
The function
$$g: \textnormal{}]0,1[^2 \to \mathbb{R}, (t,x) \mapsto \arctan(\log(e^{1/t \cdot \log^2(1/x)}+\log(e^{e^{1/t}}+2))),$$ 
is restricted log-exp-analytic.

\vs{0.3cm}
Since the global exponential function comes only locally bounded into the game one sees that a restricted log-exp-analytic function $f:\mathbb{R} \to \mathbb{R}$ fulfills the following property for all sufficiently small positive $y$. Either $f(y)$ vanishes identically or there is $c \in \mathbb{R} \setminus \{0\}$, a non-negative integer $r \in \mathbb{N}_0$ and $q_0,\ldots,q_r \in \mathbb{Q}$ such that $f(y) = c \cdot h(y) + o(h(y))$ where $h(y):=y^{q_0} \cdot (-\log(y))^{q_1} \cdot \ldots \cdot \log_{r-1}(-\log(y))^{q_r}$ (see Definition 1.13 and Proposition 3.16 in \cite{9}). A consequence is the following.

\vs{0.4cm}
{\bf Fact}

\vs{0.1cm}
Let $r \in \mathbb{R} \setminus \mathbb{Q}$. The irrational power function
$$h: \mathbb{R} \to \mathbb{R}, x \mapsto \left\{\begin{array}{ll} x^r, & x > 0, \\
0, & \textnormal{ else, } \end{array}\right.$$
is $\mathbb{R}_{\textnormal{an,exp}}$-definable, but not restricted log-exp-analytic, since one has $h \sim x^r$ as $x \searrow 0$. Because $h(x)=\exp(r\log(x))$ for $x \in \mathbb{R}_{>0}$ we see that $h|_{\mathbb{R}_{>0}}$ is restricted log-exp-analytic. 

\vs{0.3cm}
Consequently Tamm's theorem from \cite{8} is not a generalisation of the version of Tamm's theorem from \cite{3}, since in \cite{3} $\mathbb{R}_{\textnormal{an}}^{\mathbb{R}}$-definable functions are considered where $\mathbb{R}_{\textnormal{an}}^{\mathbb{R}}$ is the structure generated by all globally subanalytic functions and irrational power functions. By Miller \cite{7} the structure $\mathbb{R}_{\textnormal{an}}^{\mathbb{R}}$ is o-minimal. (See \cite{2} for the definition and properties of an o-minimal structure.) Even the structure $\mathbb{R}_{\textnormal{an,exp}}$ is o-minimal by Van den Dries \cite{1} which is a proper extension of $\mathbb{R}_{\textnormal{an}}^{\mathbb{R}}$.

\vs{0.1cm}
This article merges the results from \cite{3} and \cite{8}: We look at compositions of irrational power functions and restricted log-exp-analytic functions. Such compositions form a class of $\mathbb{R}_{\textnormal{an,exp}}$-definable functions which contains all restricted log-exp-analytic functions and all $\mathbb{R}_{\textnormal{an}}^{\mathbb{R}}$-definable functions. We call them \textit{restricted log-exp-analytic power functions}.

\vs{0.1cm}
As in \cite{8} we give differentiability results of this class of functions in the parametric setting. Thus we introduce variables $(w_1,\ldots,w_l,u_1,\ldots,u_m,z)$, where $(u_1,\ldots,u_m,z)$ is serving as the tuple of independent variables of families of functions parameterized by $w:=(w_1,\ldots,w_l)$. (The variable $z$ is needed to describe a preparation theorem for restricted log-exp-analytic power functions with respect to a single variable which is suitable for our purposes.) Then a \textit{restricted log-exp-analytic power function in $(u,z)$} on $X$ where $X \subset \mathbb{R}^l \times \mathbb{R}^m \times \mathbb{R}$ is $\mathbb{R}_{\textnormal{an,exp}}$-definable and $X_w:=\{(u,z) \in \mathbb{R}^m \times \mathbb{R} \mid (w,u,z) \in X\}$ is open for every $w \in \mathbb{R}^l$ is the composition from either side of restricted log-exp-analytic functions in $(u,z)$ and irrational power functions on $X$. A \textit{restricted log-exp-analytic function in $(u,z)$} on $X$ is the composition from either side of log-analytic functions and exponentials of locally bounded functions in $(u,z)$ ($g:X \to \mathbb{R}$ is locally bounded in $(u,z)$ if $g_w:X_w \to \mathbb{R}, (u,z) \mapsto g(w,u,z),$ is locally bounded for every $w \in \mathbb{R}^l$).

\vs{0.2cm}
One of our main goals for this article is to formulate and prove a preparation theorem for restricted log-exp-analytic power functions in $(u,z)$ (see Theorem C in \cite{9} for a precise preparation theorem for $\mathbb{R}_{\textnormal{an,exp}}$-definable functions, see $\cite{4}$ and $\cite{6}$ for original versions): a restricted log-exp-analytic power function $f:X \to \mathbb{R}, (w,u,z) \mapsto f(w,u,z),$ in $(u,z)$ where $X_w :=\{(u,z) \in \mathbb{R}^m \times \mathbb{R} \mid (w,u,z) \in X\}$ is open for $w \in \mathbb{R}^l$ can be cellwise written as \textit{$(m+1,X)$-power-restricted $(e,r)$-prepared functions} for suitable parameters $e \in \mathbb{N}_0 \cup \{-1\}$ and $r \in \mathbb{N}_0$. Here the parameter $e$ describes the maximal number of iterations of exponentials which occur in such a preparation which have the following form: each of them are exponentials of $(m+1,X)$-power-restricted $(l,r)$-prepared functions for $l<e$ which can be extended to a locally bounded function in $(u,z)$ on $X$. This information about the exponentials is described by the tuple $(m+1,X)$. The parameter $r$ describes the maximal number of iterations of the logarithm depending on $z$ which occur in every such exponential. These logarithms can be technical described by products of real powers of components of a logarithmic scale (see Definition 1.4 below for the notion of a logarithmic scale). Formally an $(m+1,X)$-power-restricted $(e,r)$-prepared function is defined as follows. 

\vspace{0.1cm}
Let $n:=l+m$, $C \subset \mathbb{R}^n \times \mathbb{R}_{\neq 0}$ be an $\mathbb{R}_{\textnormal{an,exp}}$-definable cell and let $r \in \mathbb{N}_0$. Let $t:=(w,u)$ and $\pi:\mathbb{R}^n \times \mathbb{R} \to \mathbb{R}^n, (t,z) \mapsto t,$ be the projection on the first $n$ coordinates. Let $X \subset \mathbb{R}^n \times \mathbb{R}$ be $\mathbb{R}_{\textnormal{an,exp}}$-definable with $C \subset X$ such that $X_w  \subset \mathbb{R}^{m+1}$ is open for $w \in \mathbb{R}^l$. We call a function $f:C \to \mathbb{R}, (w,u,z) \mapsto f(w,u,z),$ \textit{$(m+1,X)$-power-restricted $(-1,r)$-prepared in $z$ with center $\Theta:=(\Theta_0,\ldots,\Theta_r)$} if $f$ is the zero function. For $e \in \mathbb{N}_0$ call a function $f:C \to \mathbb{R}, (w,u,z) \mapsto f(w,u,z),$ \textit{$(m+1,X)$-power-restricted $(e,r)$-prepared in $z$ with center $\Theta$} if for $(t,z) \in C$
$$f(t,z) = \sigma a(t) \vert{y_0(t,z)}\vert^{\alpha_0} \cdot \ldots \cdot \vert{y_r(t,z)}\vert^{\alpha_r} \exp(c(t,z)) \cdot \rho(t,z)$$
where $\alpha_0,\ldots,\alpha_r \in \mathbb{R}$, $y_0=z-\Theta_0(t), y_1=\log(\vert{y_0}\vert)-\Theta_1(t),\ldots$, $\sigma \in \{-1,0,1\}$, $a=\prod_{j=1}^k h_j^{\lambda_j}$ where $k \in \mathbb{N}$, $\lambda_j \in\mathbb{R}$, $c$ can be extended to a locally bounded function in $(u,z)$ on $X$ and is itself $(m+1,X)$-power-restricted $(e-1,r)$-prepared with center $\Theta$ and $\rho(t,z)$ is a unit of a special form which we describe below. Furthermore there is $\delta>1$ such that $1/\delta<\rho<\delta$ and the functions $\Theta_0,\ldots,\Theta_r:\pi(C) \to \mathbb{R}$ and $h_j:\pi(C) \to \mathbb{R}_{>0}$ are $C$-nice functions: they are compositions of log-analytic functions and exponentials of the form $\exp(h)$ where $h$ is the component of a center of a logarithmic scale on $C$. Note that a log-analytic function on $\pi(C)$ is $C$-nice and that every $C$-nice function is definable (see \cite{9} for examples and several properties of $C$-nice functions), but the class of $C$-nice functions does not necessarily coincide with the class of definable functions: if the cell $C$ is \textit{simple}, i.e. for every $t \in \pi(C)$ there is $d_t \in \mathbb{R}_{>0} \cup \{\infty\}$ such that $C_t= \textnormal{}]0,d_t[$ (see for example Definition 2.15 in Kaiser-Opris \cite{5}), the class of $C$-nice functions coincides with the class of log-analytic ones (in \cite{5} it is shown that the center of a logarithmic scale vanishes on a simple cell).

\vspace{0.2cm}
The first goal of this paper is to prove that a restricted log-exp-analytic power function in $(u,z)$ can be indeed cellwise prepared as $(m+1,X)$-power-restricted $(e,r)$-prepared functions in $(u,z)$.

\vs{0.5cm}
{\bf Theorem A}

\vs{0.1cm}
{\it
Let $X \subset \mathbb{R}^n \times \mathbb{R}$ be $\mathbb{R}_{\textnormal{an,exp}}$-definable and let $f:X \to \mathbb{R}$ be a restricted log-exp-analytic power function in $(u,z)$. Then there are $e \in \mathbb{N}_0 \cup \{-1\}$, $r \in \mathbb{N}_0$ and an $\mathbb{R}_{\textnormal{an,exp}}$-definable cell decomposition $\mathcal{C}$ of $X_{\neq 0}$ such that for every $C \in \mathcal{C}$ there is $\Theta:=(\Theta_0,\ldots,\Theta_r)$ such that the function $f|_C$ is $(m+1,X)$-power-restricted $(e,r)$-prepared in $z$ with center $\Theta$.}

\vspace{0.1cm}
In the case of a restricted log-exp-analytic function we have a similar preparation with the difference that the function $a$ is $C$-nice and that the logarithms are rational powers of components of logarithmic scales (i.e. $\lambda_1 \ldots \lambda_k \in \mathbb{Q}$ and $\alpha_0 \ldots \alpha_r \in \mathbb{Q}$).

\vs{0.1cm}
The second goal of this paper is to give some differentiability properties for restricted log-exp-analytic power functions which are versions for Theorem A, Theorem B and Theorem C from \cite{8} for restricted log-exp-analytic power functions. These are also generalizations of the results from $\cite{3}$. 

\vs{0.5cm}
{\bf Theorem B}

\vs{0.1cm}
{\it
Let $X \subset \mathbb{R}^l \times \mathbb{R}^m$ be $\mathbb{R}_{\textnormal{an,exp}}$-definable such that $X_w$ is open for every $w \in \mathbb{R}^l$ and let $f:X \to \mathbb{R}, (w,u) \mapsto f(w,u),$ be a restricted log-exp-analytic power function in $u$. Then the following holds.
\begin{itemize}
	\item[(1)] Closedness under taking derivatives: Let $i \in \{1, \ldots , m\}$ be such that $f$ is differentiable with respect to $u_i$ on $X$. Then $\partial f/\partial u_i$ is a restricted log-exp-analytic power function in $u$.
	\item[(2)] Strong quasianalyticity: There is $N \in \mathbb{N}$ such that if $f(w,-)$ is $C^N$ for $w \in \mathbb{R}^l$ and if there is $a \in X_t$ such that all derivatives up to order $N$ vanish in $a$ then $f(w,-)$ vanishes identically.
	\item[(3)] Parametric version of Tamm's theorem: There is $M \in \mathbb{N}$ such that if $f(w,-)$ in $C^M$ at $u$ for $(w,u) \in X$ then $f(w,-)$ is real analytic at $u$. 
\end{itemize}}

\vs{0.3cm}
This paper is organised as follows. In Section 1 we pick up the most important concepts from \cite{8} and \cite {9} like log-analytic functions and their preparation theorem. In Section 2 we give a proof for Theorem A and Section 3 is devoted to the proof of Theorem B divided into three separate propositions.

\vs{0.5cm}
{\bf Notations}

\vs{0.1cm}
By $\mathbb{N}:=\{1,2,\ldots\}$ we denote the set of natural numbers and by $\mathbb{N}_0:=\{0,1,2,\ldots\}$ the set of nonnegative integers. For $m,n \in \mathbb{N}$ we denote by $\mathcal{M}(m \times n,\mathbb{R})$ respectively $\mathcal{M}(m \times n,\mathbb{Q})$ the set of $m \times n$-matrices with real respectively rational entries.

\vs{0.2cm}
For $m \in \mathbb{N}$, a set $X \subset \mathbb{R}^m$ and a set $E$ of positive real valued functions on $X$ we set $\log(E):=\{\log(g) \mid g \in E\}$.

\vs{0.2cm}
For $X \subset \mathbb{R}^n \times \mathbb{R}$ let $X_{\neq 0}:=\{(t,z) \in X \mid z \neq 0\}$. For $X \subset \mathbb{R}^l \times \mathbb{R}^m$ and $w \in \mathbb{R}^m$ we set $X_w:=\{u \in \mathbb{R}^m \mid (w,u) \in X\}$ and for a function $f:X \to \mathbb{R}, (w,u) \mapsto f(w,u),$ we set $f_w:X_w \to \mathbb{R}, u \mapsto f(w,u)$. 

\vs{0.2cm}
The reader should be familiar with basic facts about o-minimal structures from \cite{2}.

\vs{0.4cm}
{\bf Convention}

\vs{0.1cm}
Definable means $\mathbb{R}_{\textnormal{an,exp}}$-definable if not otherwise mentioned.

\section{$\mathbb{R}_{\textnormal{an,exp}}$-Definable Functions}

\subsection{Log-Analytic Functions and the Exponential \\ Number}

Compare with \cite{9}, Section 1 for a more detailed description of the content in this subsection. 

\vs{0.2cm}
Let $m \in \mathbb{N}$ and $X \subset \mathbb{R}^m$ be definable.

\vs{0.3cm}
{\bf1.1 Definition} 

\vs{0.1cm}
Let $f:X \to \mathbb{R}$ be a function.\index{log-analytic}
\begin{itemize}
	\item [(a)] Let $r \in \mathbb{N}_0$. By induction on $r$ we define that $f$ is \textbf{log-analytic of order at most} $r$.
	
	\vs{0.3cm}
	\textbf{Base case}: The function $f$ is log-analytic of order at most $0$ if there is a decomposition $\mathcal{C}$ of $X$ into finitely many definable cells such that for $C \in \mathcal{C}$ there is a globally subanalytic function $F:\mathbb{R}^m \to \mathbb{R}$ such that $f|_C = F|_C$.
	
	\vs{0.3cm}
	\textbf{Inductive step}: The function $f$ is log-analytic of order at most $r$ if the following holds: There is a decomposition $\mathcal{C}$ of $X$ into finitely many definable cells such that for $C \in \mathcal{C}$ there are $k,l \in \mathbb{N}_{0}$, a globally subanalytic function $F:\mathbb{R}^{k+l} \to \mathbb{R}$, and log-analytic functions $g_1,\ldots,g_k:C \to \mathbb{R}, h_1,\ldots,h_l:C \to \mathbb{R}_{>0}$ of order at most $r-1$ such that
	$$f|_C=F(g_1,\ldots,g_k,\log(h_1),\ldots,\log(h_l)).$$
	
	\item[(b)] Let $r \in \mathbb{N}_0$. We call $f$ \textbf{log-analytic of order}\index{log-analytic} $r$ if $f$ is log-analytic of order at most $r$ but not of order at most $r-1$.
	
	\item[(c)] We call $f$ \textbf{log-analytic}\index{log-analytic} if $f$ is log-analytic of order $r$ for some $r \in \mathbb{N}_0$.
\end{itemize}

\vs{0.3cm}
{\bf1.2 Definition}

\vs{0.1cm}
Let $f:X \to \mathbb{R}$ be a function. Let $E$ be a set of positive definable functions on $X$.
\begin{itemize}
	\item [(a)] By induction on $e \in \mathbb{N}_0$ we define that $f$ has \textbf{exponential number at most $e$ with respect to $E$}\index{exponential!number}.
	
	\vs{0.2cm}
	{\bf Base Case}: The function $f$ has exponential number at most $0$ with respect to $E$ if $f$ is log-analytic.
	
	\vs{0.2cm}
	{\bf Inductive Step}: The function $f$ has exponential number at most $e$ with respect to $E$ if the following holds: There are $k,l \in \mathbb{N}_0$, functions $g_1,\ldots,g_k:X \to \mathbb{R}$ and $h_1,\ldots,h_l:X \to \mathbb{R}$ with exponential number at most $e-1$ with respect to $E$ and a log-analytic function $F:\mathbb{R}^{k+l} \to \mathbb{R}$ such that
	$$f=F(g_1,\ldots,g_k,\exp(h_1),\ldots,\exp(h_l))$$
	and $\exp(h_1),\ldots,\exp(h_l) \in E$.
	
	\item [(b)] Let $e \in \mathbb{N}_0$. We say that $f$ has \textbf{exponential number $e$ with respect to $E$}\index{exponential!number} if $f$ has exponential number at most $e$ with respect to $E$ but not at most $e-1$ with respect to $E$.
	
	\item [(c)] We say that $f$ \textbf{can be constructed from $E$}\index{constructed from} if there is $e \in \mathbb{N}_0$ such that $f$ has exponential number $e$ with respect to $E$. 
\end{itemize}

\vs{0.5cm}
{\bf1.3 Remark}

\vs{0.1cm}
Let $e \in \mathbb{N}_0$. Let $E$ be a set of positive definable functions on $X$.
\begin{itemize}
	\item[(1)] Let $f:X \to \mathbb{R}$ be a function with exponential number at most $e$ with respect to $E$. Then $\exp(f)$ has exponential number at most $e+1$ with respect to $E \cup \{\exp(f)\}$.
	\item[(2)] Let $s \in \mathbb{N}_0$. Let $f_1,\ldots,f_s:X \to \mathbb{R}$ be functions with exponential number at most $e$ with respect to $E$ and let $F:\mathbb{R}^s \to \mathbb{R}$ be log-analytic. Then $F(f_1,\ldots,f_s)$ has exponential number at most $e$ with respect to $E$. 
\end{itemize}

\subsection{A Preparation Theorem for Log-Analytic Functions}

Compare with \cite{9}, Section 2 for a more detailed description of the content in this subsection.

\vs{0.2cm}
Let $n \in \mathbb{N}$. Let $t$ range over $\mathbb{R}^n$ and $z$ over $\mathbb{R}$. We fix a definable set $C \subset \mathbb{R}^n \times \mathbb{R}$.

\vs{0.5cm}
{\bf1.4 Definition} (\cite{9} Section 2.1)

\vs{0.1cm}
Let $r \in \mathbb{N}_0$. A tuple $\mathcal{Y}:=(y_0,\ldots,y_r)$ of functions on $C$ is called an \textbf{$r$-logarithmic scale} on $C$ with \textbf{center} $\Theta=(\Theta_0,\ldots,\Theta_r)$ if the following holds:
\begin{itemize}
	\item[(a)] $y_j>0$ or $y_j<0$ for every $j \in \{0,\ldots,r\}$.
	\item[(b)] $\Theta_j$ is a definable function on $\pi(C)$ for every $j \in \{0,\ldots,r\}$.
	\item[(c)] We have $y_0(t,z)=z-\Theta_0(t)$ and inductively $y_j(t,z)=\log(\vert{y_{j-1}(t,z)}\vert) - \Theta_j(t)$ for every $j \in \{1,\ldots,r\}$ and all $(t,z) \in C$.
	\item[(d)] Either there is $\epsilon_0 \in \textnormal{}]0,1[$ such that $0<\vert{y_0(t,z)}\vert < \epsilon_0\vert{z}\vert$ for all $(t,z) \in C$ or $\Theta_0=0$, and for every $j \in \{1,\ldots,r\}$ either there is $\epsilon_j \in \textnormal{}]0,1[$ such that $0<\vert{y_j(t,z)}\vert<\epsilon_j\vert{\log(\vert{y_{j-1}(t,z)}\vert)}\vert$ for all $(t,z) \in C$ or $\Theta_j=0$.
\end{itemize}

\vs{0.2cm}
For a logarithmic scale $(y_0,\ldots,y_r)$ on a definable set $C$ and $\alpha \in \mathbb{R}^{r+1}$ we often write $\vert{\mathcal{Y}(t,z)}\vert^{\otimes \alpha}$ instead of $\prod_{j=0}^r\vert{y_j(t,z)}\vert^{\alpha_j}$ where $(t,z) \in C$.

\vs{0.4cm}
{\bf1.5 Definition} (\cite{9} Section 2.3)

\vs{0.1cm}
We call $g:\pi(C) \to \mathbb{R}$ a \textbf{$C$-heir} if there is $l \in \mathbb{N}_0$, an $l$-logarithmic scale $\hat{\mathcal{Y}}$ with center $(\hat{\Theta}_0,\ldots,\hat{\Theta}_l)$ on $C$, and $j \in \{1,\ldots,l\}$ such that $g=\exp(\Theta_j)$.

\vs{0.5cm}
{\bf1.6 Definition} (\cite{9} Section 2.3)

\vs{0.1cm}
We call $g:\pi(C) \to \mathbb{R}$ \textbf{$C$-nice} if there is a set $E$ of $C$-heirs such that $g$ can be constructed from $E$.

\vs{0.5cm}
Note that the class of log-analytic functions on $\pi(C)$ can be a proper subclass of the class of $C$-nice functions (compare with Example 2.39 in \cite{9}). In the following we give the definition from $\cite{9}$ for log-analytically prepared functions with the difference that we also allow real exponents for the iterations of the logarithms. This is needed to describe preparations of restricted log-exp-analytic power functions on simple cells in an effective way.

\vs{0.5cm}
{\bf1.7 Definition}

\vs{0.1cm}
Let $r \in \mathbb{N}_0$. Let $g:C \to \mathbb{R}$ be a function. We say that $g$ is \textbf{$r$-real-log-analytically prepared in $z$ with center $\Theta$} if
$$g(t,z)=a(t) \vert{\mathcal{Y}(t,z)}\vert^{\otimes \alpha}\rho(t,z)$$
for all $(t,z) \in C$ where $a$ is a definable function on $\pi(C)$ which vanishes identically or has no zero, $\mathcal{Y}=(y_0,\ldots,y_r)$ is an $r$-logarithmic scale with center $\Theta$ on $C$, $\alpha \in \mathbb{R}^{r+1}$ and the following holds for $\rho$. There is $s \in \mathbb{N}$ such that $\rho=v \circ \phi$ where $v$ is a power series which converges on an open neighbourhood of $[-1,1]^s$ with $v([-1,1]^s) \subset \mathbb{R}_{>0}$ and $\phi:=(\phi_1,\ldots,\phi_s):C \to [-1,1]^s$ is a function of the form 
$$\phi_j(t,z):=b_j(t)\vert{\mathcal{Y}(t,z)}\vert^{\otimes \gamma_j}$$
for $j \in \{1,\ldots,s\}$ and $(t,z) \in C$ where $b_j:\pi(C) \to \mathbb{R}$ is definable for $j \in \{1,\ldots,s\}$ and $\gamma_j:=(\gamma_{j0},\ldots,\gamma_{jr}) \in \mathbb{R}^{r+1}$. 
We call $a$ \textbf{coefficient} and $b:=(b_1,\ldots,b_s)$ a tuple of \textbf{base functions} for $f$. An \textbf{LA-preparing tuple} for $f$ is then
$$\mathcal{J}:=(r,\mathcal{Y},a,\alpha,s,v,b,\Gamma)$$
where
$$\Gamma:=\left(\begin{array}{cccc}
\gamma_{10}&\cdot&\cdot&\gamma_{1r}\\
\cdot&& &\cdot\\
\cdot&& &\cdot\\
\gamma_{s0}&\cdot&\cdot&\gamma_{sr}\\
\end{array}\right)\in \mathcal{M}\big(s\times (r+1),\mathbb{R}).$$

If $\alpha,\gamma_1,\ldots,\gamma_s \in \mathbb{Q}^{r+1}$ we say that $g$ is \textbf{$r$-log-analytically prepared in $z$ with center $\Theta$}.

\vs{0.5cm}
The following preparation theorem for log-analytic functions has been established in \cite{9}.

\vs{0.5cm}
{\bf1.8 Fact} (\cite{9} Theorem A)

\vs{0.1cm}
{\it
	Let $m \in \mathbb{N}$, $r \in \mathbb{N}_0$. Let $X \subset \mathbb{R}^n \times \mathbb{R}$ be definable. Let $f_1,\ldots,f_m:X \to \mathbb{R}$ be log-analytic functions of order at most $r$. Then there is a definable cell decomposition $\mathcal{C}$ of $X_{\neq 0}$ such that $f_1|_C,\ldots,f_m|_C$ are $r$-log-analytically prepared in $z$ with $C$-nice coefficient, $C$-nice base functions and common $C$-nice center for $C \in \mathcal{C}$.}

\section{Restricted Log-Exp-Analytic \\
	 Power Functions}

\subsection{Basic Facts and Definitions}

The main results of this paper are formulated in the parametric setting. So we set up the concept of restricted log-exp-analytic power functions in single variables.

\vs{0.3cm}
Let $l,m \in \mathbb{N}_0$. Let $w$ range over $\mathbb{R}^l$ and $u$ over $\mathbb{R}^m$. We fix definable sets $C,X \subset \mathbb{R}^l \times \mathbb{R}^m$ with $C \subset X$. Suppose that $X_w$ is open for every $w \in \mathbb{R}^l$. Let $\pi_l:\mathbb{R}^l \times \mathbb{R}^m \to \mathbb{R}^l, (w,u) \mapsto w$.

\vs{0.5cm}
{\bf2.1 Definition}

\vs{0.1cm}
We call a function $g:\mathbb{R} \to \mathbb{R}$ \textbf{power function} if there is $\chi \in \mathbb{R}$ such that $g(x)=x^\chi$ for $x \in \mathbb{R}_{>0}$ and $g(x)=0$ otherwise.

\vspace{0.3cm}
Note that power functions are definable since $x^\chi=e^{\chi \log(x)}$ for every $c \in \mathbb{R}$ and $x \in \mathbb{R}_{>0}$. Our next aim is to define restricted log-exp-analytic power functions formally which are compositions of log-analytic functions, exponentials of locally bounded functions and power functions. In the sense of Definition 1.2(c) they are precisely those functions which can be constructed from a set $E$ of positive definable functions such that every $g \in \log(E)$ is locally bounded or $g=\chi \log(h)$ for a constant $\chi \in \mathbb{R}$ and a positive function $h$ which can also be constructed from $E$. For convenience we call such a set $E$ a \textit{LoPo-set}.

\vs{0.5cm}
{\bf2.2 Definition}

\vs{0.1cm}
Let $E$ be a set of positive definable functions on $C$. We call $E$ a \textbf{LoPo-set on $C$ in $u$ with reference set $X$} if the following holds: Let $e \in \mathbb{N}_0$ and let $g \in \log(E)$ be with exponential number at most $e$ with respect to $E$. Then $g$ is locally bounded in $u$ with reference set $X$ (i.e. there is a definable function $\tilde{g}:X \to \mathbb{R}$ with $\tilde{g}|_C=g$ where $\tilde{g}_w$ is locally bounded for $w \in \pi_l(X)$) or there is a function $h:C \to \mathbb{R}_{>0}$ which has exponential number at most $e$ with respect to $E$ and a constant $\chi \in \mathbb{R}$ such that $g=\chi \log(h)$.

\vs{0.3cm}
{\bf2.3 Remark}

\vs{0.1cm}
Let $E$ be a set of positive definable functions on $C$. Let $Y \subset \mathbb{R}^l \times \mathbb{R}^m$ be definable with $X \subset Y$ such that $Y_w$ is open for every $w \in \mathbb{R}^l$. Let $E$ be a LoPo-set in $u$ with reference set $Y$. Then $E$ is a LoPo-set in $u$ with reference set $X$.

\vs{0.3cm}
{\bf Proof}

\vs{0.1cm}
This follows from the following fact. Let $g:C \to \mathbb{R}$ be locally bounded in $u$ with reference set $Y$. Then $g:C \to \mathbb{R}$ is locally bounded in $u$ with reference set $X$. \hfill$\blacksquare$ 

\vs{0.5cm}
{\bf2.4 Definition}

\vs{0.1cm}
Let $f:C \to \mathbb{R}$ be a function.
\begin{itemize}
	\item [(a)]
	Let $e \in \mathbb{N}_0$. We say that $f$ is a \textbf{restricted log-exp-analytic power function (restricted log-exp-analytic function) in $u$ of order (at most) $e$ with reference set $X$} if $f$ has exponential number (at most) $e$ with respect to a LoPo-set $E$ in $u$ (with respect to a set $E$ of exponentials of locally bounded functions in $u$) with reference set $X$ on $C$.  
	
	\item [(b)]
	We say that $f$ is a \textbf{restricted log-exp-analytic power function (restricted log-exp-analytic function) in $u$ with reference set $X$} if $f$ can be constructed from a LoPo-set $E$ in $u$ (from a set $E$ of exponentials of locally bounded functions in $u$) with reference set $X$ on $C$, i.e. there is $e \in \mathbb{N}_0$ and a LoPo-set $E$ in $u$ (a set $E$ of exponentials of locally bounded functions in $u$) on $C$ with reference set $X$ such that $f$ has exponential number (at most) $e$ with respect to $E$.
\end{itemize}

\vs{0.3cm}
{\bf2.5 Remark}
\begin{itemize}
	\item [(1)]
	The log-analytic functions are precisely the restricted log-exp-analytic power functions in $u$ of order (at most) $0$.
	
	\item [(2)]
	A restricted log-exp-analytic function $f:C \to \mathbb{R}$ in $u$ with reference set $X$ is a restricted log-exp-analytic power function in $u$ with reference set $X$.
\end{itemize}

\vs{0.3cm}
{\bf2.6 Example}

\vs{0.1cm} 
Let $\chi \in \mathbb{R} \setminus \mathbb{Q}$. The irrational power function
$$f: \mathbb{R} \to \mathbb{R}, u \mapsto \left\{\begin{array}{ll} u^\chi, & u > 0, \\
0, & \textnormal{ else, } \end{array}\right.$$
is a restricted log-exp-analytic power function (of order (at most) $1$) in $u$ with reference set $\mathbb{R}$.

\vs{0.4cm}
{\bf Proof}

\vs{0.1cm}
This is immediately seen with the fact that $f(u)=\exp(\chi \log(u))$ for $u \in \mathbb{R}_{>0}$ and $f(u)=0$ otherwise: let 
\[
g:\mathbb{R} \to \mathbb{R}, u \mapsto \left\{\begin{array}{ll} \exp(\chi \log(u)), & u > 0, \\
	1, & \text{else,} \end{array}\right.
\]
and let $E:=\{g\}$. Then $E$ is a LoPo-set in $u$ with reference set $\mathbb{R}$, since $\log(g)=\chi \log(h)$ for the log-analytic function $h:\mathbb{R} \to \mathbb{R}_{>0}$ with $h(u)=u$ for $u>0$ and $h=1$ otherwise. Let
\[
G:\mathbb{R}^2 \to \mathbb{R}, (x_1,x_2) \mapsto \left\{\begin{array}{ll} x_1, & x_2 > 0, \\
	0, & \text{else.} \end{array}\right.
\]
Then $G$ is log-analytic (and even globally subanalytic). Since $f(u)=G(g(u),u)$ for $u \in \mathbb{R}$ we see that $f$ is a restricted log-exp-analytic power function of order $1$ in $u$ with reference set $\mathbb{R}$ (since $f$ is not log-analytic).
\hfill$\blacksquare$

\vs{0.5cm}
{\bf2.7 Remark}

\vs{0.1cm}
Let $e \in \mathbb{N}_0$. Let $Y \subset \mathbb{R}^l \times \mathbb{R}^m$ be definable with $X \subset Y$ such that $Y_w$ is open for every $w \in \mathbb{R}^l$. Let $f:C \to \mathbb{R}$ be a restricted log-exp-analytic power function in $u$ of order at most $e$ with reference set $Y$. Then $f$ is a restricted log-exp-analytic power function in $u$ of order at most $e$ with reference set $X$.

\vs{0.3cm}
{\bf Proof}

\vs{0.1cm}
This is directly seen with Remark 2.3. \hfill$\blacksquare$

\vs{0.5cm}
{\bf2.8 Remark}

\vs{0.1cm}
Let $k \in \mathbb{N}$. For $j \in \{1,\ldots,k\}$ let $f_j:C \to \mathbb{R}$ be a restricted log-exp-analytic power function in $u$ with reference set $X$. Let $F:\mathbb{R}^k \to \mathbb{R}$ be log-analytic. Then 
$$C \to \mathbb{R}, u \mapsto F(f_1(u),\ldots,f_k(u)),$$
is a restricted log-exp-analytic power function in $u$ with reference set $X$.

\vs{0.3cm}
{\bf Proof}

\vs{0.1cm}
Note that $f_j$ can be constructed from a set $E_j$ of positive definable functions which is a LoPo-set in $u$ with reference set $X$ for $j \in \{1,\ldots,k\}$. Then $E:=E_1 \cup \ldots \cup E_k$ is a LoPo-set in $u$ with reference set $X$ and for $j \in \{1,\ldots,k\}$ the function $f_j$ can be constructed from $E$. With Remark 1.3(2) we are done. \hfill$\blacksquare$

\vs{0.5cm}
{\bf2.9 Remark}

\vs{0.1cm}
Let $C_1,C_2 \subset \mathbb{R}^m$ be disjoint and definable with $C_1 \cup C_2 = C$. For $j \in \{1,2\}$ let $f_j:C_j \to \mathbb{R}$ be a restricted log-exp-analytic power function in $u$ with reference set $X$. Then
$$f:C \to \mathbb{R}, (w,u) \mapsto \left\{\begin{array}{ll} f_1(w,u) , & (w,u) \in C_1,  \\
f_2(w,u), & (w,u) \in C_2, \end{array}\right.$$
is a restricted log-exp-analytic power function in $u$ with reference set $X$.

\vspace{0.2cm}
{\bf Proof}

\vspace{0.1cm}
For $j \in \{1,2\}$ let $E_j$ be a LoPo-set on $C_j$ in $x$ with reference set $X$ such that $f_j$ can be constructed from $E_j$. For $j \in \{1,2\}$ let
$$\tilde{E}_j:=\{\delta:C \to \mathbb{R} \mid \textnormal{$\delta$ is a function with } \delta|_{C_j} \in E_j \textnormal{ and } 1 \text{ otherwise}\}.$$
Then $\tilde{E}_j$ is a LoPo-set on $C$ in $u$ with reference set $X$: let $g \in \log(\tilde{E}_j)$ be with exponential number at most $e$ with respect to $\tilde{E}_j$. Then $g|_{C_j} \in \log(E_j)$ has exponential number at most $e$ with respect to $\tilde{E}_j|_{C_j}=E_j$. 

If $g|_{C_j}$ is of the form $\chi \log(h)$ then $g$ is of the form $\chi \log(\tilde{h})$ with $\tilde{h}(w,u)=h(w,u)$ for $(w,u) \in C_j$ and $\tilde{h}(w,u)=1$ otherwise. Note that $\tilde{h}$ has exponential number at most $e$ with respect to $\tilde{E}_j$. 

If $g$ is locally bounded in $u$ with reference set $X$ then $\tilde{g}:C \to \mathbb{R}$ with $\tilde{g}(w,u)=g(w,u)$ for $(w,u) \in C_j$ and $0$ otherwise is also locally bounded in $u$ with reference set $X$. Therefore $E:=\tilde{E}_1 \cup \tilde{E}_2$ is a LoPo-set in $u$ with reference set $X$ from which $f$ can be constructed. \hfill$\blacksquare$

\vs{0.5cm}
{\bf2.10 Definition}

\vs{0.1cm}
A function $f:X \to \mathbb{R}$ is called a \textbf{restricted log-exp-analytic power function in $u$} if $f$ is a restricted log-exp-analytic power function in $u$ with reference set $X$.

\vs{0.5cm}
{\bf2.11 Remark}

\vs{0.1cm}
Let $k \in \mathbb{N}_0$. Let $v:=(v_1,\ldots,v_k)$ range over $\mathbb{R}^k$. Let $g:\mathbb{R}^k \to \mathbb{R}^m$ be log-analytic and continuous. Let
$$V:=\{(w,u,v) \in X \times \mathbb{R}^k \mid (w,u+g(v)) \in X\}.$$
Let $f:X \to \mathbb{R}, (w,u) \mapsto f(w,u),$ be a restricted log-exp-analytic power function in $u$. 
Then $F:V \to \mathbb{R}, (w,u,v) \mapsto f(w,u+g(v)),$ is a restricted log-exp-analytic power function in $(u,v)$. 

\vs{0.4cm}
{\bf Proof}

\vs{0.1cm}
Note that $V_w$ is open in $\mathbb{R}^m \times \mathbb{R}^k$ for every $w \in \mathbb{R}^l$. Let $E$ be a LoPo-set in $u$ with reference set $X$ such that $f$ can be constructed from $E$. Consider
$$\tilde{E}:=\{V \to \mathbb{R}_{>0}, (w,u,v) \mapsto h(w,u+g(v)) \mid h \in E\}.$$  
Note that $F$ can be constructed from $\tilde{E}$. We show that $\tilde{E}$ is a LoPo-set in $(u,v)$ with reference set $V$ and we are done. Let $\beta \in \tilde{E}$. Then there is $h \in E$ with $\beta(w,u,v)=h(w,u+g(v))$ for $(w,u,v) \in V$. Let $e \in \mathbb{N}$ be such that $h$ has exponential number at most $e$ with respect to $E$.

\vs{0.2cm}
{\bf Case 1}: Let $h$ be locally bounded in $u$ with reference set $X$. Then $\beta$ is locally bounded in $(u,v)$ with reference set $V$ by the claim in the proof of Remark 2.10 in \cite{8}.

\vs{0.2cm}
{\bf Case 2}: Let $\chi \in \mathbb{R}$ be a constant and $\eta:X \to \mathbb{R}_{>0}$ be with exponential number at most $e-1$ with respect to $E$ such that $h=\exp(\chi \log(\eta))$. We obtain $$\beta(w,u,v)=h(w,u+g(v))=\exp(\chi \log(\eta(w,u+g(v))))$$
for $(w,u,v) \in V$. Note that $V \to \mathbb{R}_{>0}, (w,u,v) \mapsto \eta(w,u+g(v)),$ has exponential number at most $(e-1)$ with respect to $\tilde{E}$. This finishes the proof. \hfill$\blacksquare$

\subsection{A Preparation Theorem for Restricted \\
	Log-Exp-Analytic Power Functions}

In this section we give a preparation theorem for restricted log-exp-analytic power functions. Our considerations start with Theorem $B$ from \cite{9}. 

\vs{0.2cm}
Let $m,l \in \mathbb{N}_0$. Let $w$ range over $\mathbb{R}^l$ and $u$ over $\mathbb{R}^m$. Here $(u_1,\ldots,u_m,z)$ is serving as the tuple of independent variables of families of functions parameterized by $w:=(w_1,\ldots,w_l)$. Furthermore we fix definable sets $C,X \subset \mathbb{R}^l \times \mathbb{R}^m$ with $C \subset X$ such that $X_w$ is open for $w \in \mathbb{R}^n$.

\vs{0.5cm}
{\bf2.12 Definition}

\vs{0.1cm}
Let $f:C \to \mathbb{R}$ be definable. Suppose that $f(x)>0$ for every $x \in C$, $f(x)<0$ for every $x \in C$ or $f=0$. Then $f$ is \textbf{a finite product of powers of definable functions $g_1,\ldots,g_k:C \to \mathbb{R}_{>0}$} for $k \in \mathbb{N}$ if there are $\lambda_1,\ldots,\lambda_k \in \mathbb{R}$ and $\sigma \in \{-1,0,1\}$ such that $f=\sigma \prod_{j=1}^k g_j^{\lambda_j}$.

\vs{0.5cm}
{\bf2.13 Definition}

\vs{0.1cm}
Let $f:C \to \mathbb{R}, (w,u) \mapsto f(w,u),$ be a function. By induction on $e \in \mathbb{N}_0 \cup \{-1\}$ we define that $f$ is \textbf{$(m,X)$-power-restricted $e$-prepared}. To this preparation we associate a \textbf{finite set of log-analytic functions $L$ on $C$ which ''occur'' in this preparation}. 

\vs{0.2cm}
$e=-1$: The function $f$ is $(m,X)$-power-restricted $(-1)$-prepared if $f$ is the zero function. Then $L:=\{0\}$.

\vs{0.2cm}
$e-1 \to e$: The function $f$ is $(m,X)$-power-restricted $e$-prepared if the following holds. There is $s \in \mathbb{N}$ such that
$$f(w,u)=a(w,u)\exp(c(w,u))v(b_1(w,u)\exp(d_1(w,u)),\ldots,b_s(w,u)\exp(d_s(w,u)))$$
for $(w,u) \in C$ where $a,b_1,\ldots,b_s$ are finite products of powers of log-analytic functions, $c,d_1,\ldots,d_s$ are locally bounded in $u$ with reference set $X$ and are $(m,X)$-power-restricted $(e-1)$-prepared. Additionally we have $b_j(w,u)\exp(d_j(w,u) \in [-1,1]$ for $(w,u) \in C$ and $v$ is a power series which converges on an open neighbourhood of $[-1,1]^s$ with $v([-1,1]^s) \subset \mathbb{R}_{>0}$. Suppose that for $c$ and $d_1,\ldots,d_s$ corresponding sets of log-analytic functions $L_c,L_{d_1},\ldots,L_{d_s}$ have already been defined. Let $b_0:=a$. For $j \in \{0,\ldots,s\}$ let $\sigma_j \in \{-1,0,1\}$ and $\lambda_{j0},\ldots,\lambda_{jk} \in \mathbb{R}$, $h_{j0},\ldots,h_{jk}:C \to \mathbb{R}_{>0}$ be log-analytic with
$$b_j= \sigma_j \prod_{i=1}^k h_{ji}^{\lambda_{ji}}$$
where $k \in \mathbb{N}$. We set 
$$L:=L_c \cup L_{d_1} \cup \ldots \cup L_{d_s} \cup \{h_{ji} \mid j \in \{0,\ldots,s\}, i \in \{1,\ldots,k\}\}.$$

\vs{0.4cm}
\textbf{Convention}

\vs{0.1cm}
For a set $E$ of positive definable functions on $X$ we say that $f:X \to \mathbb{R}$ has exponential number at most $-1$ with respect to $E$ if $f$ is the zero function.

\vs{0.5cm}
{\bf2.14 Proposition}

\vs{0.1cm}
{\it Let $e \in \mathbb{N}_0$. Let $f:X \to \mathbb{R}$ be a restricted log-exp-analytic power function in $u$ of order at most $e$. Then there is a decomposition $\mathcal{C}$ of $X$ into finitely many definable cells such that for every $C \in \mathcal{C}$ the function $f|_C$ is $(m,X)$-power-restricted $e$-prepared.}

\vs{0.3cm}
{\bf Proof}

\vs{0.1cm}
Let $e \in \mathbb{N}_0 \cup \{-1\}$ and $E$ be a LoPo-set in $u$ with reference set $X$ such that $f$ has exponential number at most $e$ with respect to $E$. We proceed by induction on $e$. For $e=-1$ the assertion is clear. 

\vs{0.1cm}
$e-1 \to e$: There is a decomposition $\mathcal{D}$ of $X$ into finitely many definable cells such that for every $D \in \mathcal{D}$ there is $s \in \mathbb{N}$ such that
$$f(w,u)=a(w,u)\exp(d_0(w,u))v(b_1(w,u)\exp(d_1(w,u)),\ldots,b_s(w,u)\exp(d_s(w,u)))$$
for $(w,u) \in D$ where $a,b_1,\ldots,b_s:D \to \mathbb{R}$ are log-analytic and $d_0,d_1,\ldots,d_s:D \to \mathbb{R}$ are finite $\mathbb{Q}$-linear combinations of functions from $\log(E)$ which have exponential number at most $(e-1)$ with respect to $E$. Additionally $b_j(w,u)\exp(d_j(w,u)) \in [-1,1]$ for $(w,u) \in D$ and $v$ is a power series which converges absolutely on an open neighbourhood of $[-1,1]^s$ with $v([-1,1]^s) \subset \mathbb{R}_{>0}$ (see Theorem B in \cite{9}). Fix $D \in \mathcal{D}$ with the corresponding preparation for $f|_D$. Note that there are locally bounded $\delta_0,\ldots,\delta_s:D \to \mathbb{R}$ in $u$ with reference set $X$ which have exponential number at most $(e-1)$ with respect to $E$ and functions $\zeta_0,\ldots,\zeta_s:D \to \mathbb{R}$ such that $d_j=\delta_j+\zeta_j$  for $j \in \{0,\ldots,s\}$ and the following holds. There are $k \in \mathbb{N}$ and constants $\chi_{j1},\ldots,\chi_{jk} \in \mathbb{R}$ and positive functions $\eta_{j1},\ldots,\eta_{jk}:D \to \mathbb{R}_{>0}$ which have exponential number at most $e-1$ with respect to $E$ such that $\zeta_j = \sum_{i=1}^k \chi_{ji} \log(\eta_{ji})$. So we obtain
$$f|_D=a(\prod_{i=1}^k\eta_{0i}^{\chi_{0i}}) \exp(\delta_0)v(b_1(\prod_{i=1}^k\eta_{1i}^{\chi_{1i}})\exp(\delta_1),\ldots,b_s(\prod_{i=1}^k\eta_{si}^{\chi_{si}})\exp(\delta_s)).$$
Now we use the inductive hypothesis on $\eta_{ji}$ and find a decomposition $\mathcal{A}$ of $D$ into finitely many definable cells such that for $j \in \{0,\ldots,s\}$ and $i \in \{1,\ldots,k\}$ we have that $\eta_{ji}$ is $(m,X)$-power-restricted $(e-1)$-prepared, i.e. for $(w,u)  \in A$ 
$$\eta_{ji}(w,u)=\hat{a}_{ji}(w,u)\exp(\nu_{j,i,0}(w,u)) \cdot $$ $$\hat{v}_{ji}(\hat{b}_{j,i,1}(w,u)\exp(\nu_{j,i,1}(w,u)),\ldots,\hat{b}_{j,i,\hat{s}}(w,u)\exp(\nu_{j,i,\hat{s}}(w,u)))$$

\vs{0.2cm}
where $\nu_{j,i,0},\ldots,\nu_{j,i,\hat{s}}:A \to \mathbb{R}$ are locally bounded in $u$ with reference set $X$, the functions $\hat{a}_{ji},\hat{b}_{j,i,1},\ldots,\hat{b}_{j,i,\hat{s}}:A \to \mathbb{R}$ are finite products of powers of log-analytic functions and $\hat{v}_{ji}$ is a power series which converges absolutely on an open neighbourhood of $[-1,1]^{\hat{s}}$ with $\hat{v}_{ji}([-1,1]^{\hat{s}}) \subset \mathbb{R}_{>0}$. (By redefining the single $v_{ij}$ we may assume that $\hat{s}$ does not depend on $j$.) Note also that $\hat{a}_{ji}$ is positive.

\vs{0.2cm}
Fix $A \in \mathcal{A}$ and the corresponding preparation for $\eta_{ji}|_A$. Let $$\beta_{ji}:=\hat{v}_{ji}(\hat{b}_{j,i,1}\exp(\nu_{j,i,1}),\ldots,\hat{b}_{j,i,\hat{s}}\exp(\nu_{j,i,\hat{s}})).$$
For $j \in \{0,\ldots,s\}$ let 
$$\omega_j:=\prod_{i=1}^k \beta_{ji}^{\chi_{ji}}, \textnormal{  } \kappa_j:=\sum_{i=1}^k \chi_{ji}\nu_{j,i,0}, \textnormal{  } \mu_j := \prod_{i=1}^{k}\hat{a}_{ji}^{\chi_{ji}}.$$
Note that $\hat{v}_{ji}^{\chi_{ji}}$ is a power series which converges absolutely on an open neighbourhood of $[-1,1]^s$ with $\hat{v}_{ji}^{\chi_{ji}}([-1,1]^s) \subset \mathbb{R}_{>0}$ (by using the exponential series, the logarithmic series and the fact that $\hat{v}_{ji}^{\chi_{ji}} = \exp(\chi_{ji}\log(\hat{v}_{ji}))$). We obtain
$$f|_A=a\mu_0e^{\delta_0+\kappa_0}\omega_0v(b_1\mu_1e^{\delta_1+\kappa_1}\omega_1,\ldots,b_s\mu_se^{\delta_s+\kappa_s}\omega_s).$$
Note that $a\mu_0$ and $b_j\mu_j$ for $j \in \{1,\ldots,s\}$ are finite product of powers of log-analytic functions. Additionally $\delta_j+\kappa_j$ is locally bounded in $u$ with reference set $X$ and has exponential number at most $e-1$ with respect to $E$ for $j \in \{0,\ldots,s\}$. So by the inductive hypothesis on $\delta_j+\kappa_j$ for $j \in \{0,\ldots,\hat{s}\}$ we find a decomposition $\mathcal{B}$ of $A$ into finitely many definable cells such that for every $B \in \mathcal{B}$ we have that $(\delta_j+\kappa_j)|_B$ is $(m,X)$-power-restricted $(e-1)$-prepared for $j \in \{0,\ldots,s\}$. We are done by composition of power series. \hfill$\blacksquare$

\vs{0.5cm}
For the rest of Section 2.2 let $n \in \mathbb{N}_0$ be with $n=l+m$. Let $z$ range over $\mathbb{R}$. Now $(u_1,\ldots,u_m,z)$ is serving as the tuple of independent variables of families of functions parameterized by $w:=(w_1,\ldots,w_l)$. Let $t:=(w,u)$ and let $\pi:\mathbb{R}^{n+1} \to \mathbb{R}^n, (t,z) \mapsto t$. Now let $C,X \subset \mathbb{R}^n \times \mathbb{R}$ be definable sets with $C \subset X$ such that $X_w$ is open for every $w \in \mathbb{R}^l$. 

\vs{0.5cm}
{\bf2.15 Definition}

\vs{0.1cm}
Let $e \in \mathbb{N}_0 \cup \{-1\}$ and $r \in \mathbb{N}_0$. By induction on $e \in \mathbb{N}_0 \cup \{-1\}$ we define that $f:C \to \mathbb{R}, (w,u,z) \mapsto f(w,u,z),$ is \textbf{$(m+1,X)$-power-restricted $(e,r)$-prepared in $z$} and associate a \textbf{preparing tuple} to this preparation. 

\vs{0.2cm}
$e=-1$: We call $f$ $(m+1,X)$-power-restricted $(-1,r)$-prepared in $z$ if $f$ is the zero function. A preparing tuple for $f$ is then $(0)$.

\vs{0.2cm}
$e-1 \to e$: We call $f$ $(m+1,X)$-power-restricted $(e,r)$-prepared in $z$ if for $(t,z) \in C$
$$f(t,z)=a(t)\vert{\mathcal{Y}(t,z)}\vert^{\otimes \alpha} \exp(d_0(t,z))\rho(t,z)$$
where $a:\pi(C) \to \mathbb{R}$ is a finite product of powers of $C$-nice functions, $\alpha \in \mathbb{R}^{r+1}$, $d_0:C \to \mathbb{R}$ is locally bounded in $(u,z)$ with reference set $X$ and is $(m+1,X)$-power-restricted $(e-1,r)$-prepared in $z$. Additionally $\rho:C \to \mathbb{R}$ is of the following form. There is $s \in \mathbb{N}$ such that $\rho=v \circ \phi$ where $\phi:=(\phi_1,\ldots,\phi_s):C \to [-1,1]^s$ with 
$$\phi_j(t,z)=b_j(t)\vert{\mathcal{Y}(t,z)}\vert^{\otimes \gamma_j}\exp(d_j(t,z))$$
for $j \in \{1,\ldots,s\}$ where $\gamma_j \in \mathbb{R}^{r+1}$, $b_j:\pi(C) \to \mathbb{R}$ is a finite product of powers of $C$-nice functions, $d_j:C \to \mathbb{R}$ is $(m+1,X)$-power-restricted $(e-1,r)$-prepared in $z$ and locally bounded in $(u,z)$ with reference set $X$ and $v$ is a power series which converges absolutely on an open neighbourhood of $[-1,1]^s$ with $v([-1,1]^s) \subset \mathbb{R}_{>0}$. A preparing tuple for $f$ is then 
$$(r,\mathcal{Y},a,\exp(d_0),\alpha,s,v,b,\exp(d),\Gamma)$$ 
with $b:=(b_1,\ldots,b_s)$, $\exp(d):=(\exp(d_1),\ldots,\exp(d_s))$ and
$$\Gamma:=\left(\begin{array}{cccc}
\gamma_{10}&\cdot&\cdot&\gamma_{1r}\\
\cdot&& &\cdot\\
\cdot&& &\cdot\\
\gamma_{s0}&\cdot&\cdot&\gamma_{sr}\\
\end{array}\right)\in \mathcal{M}\big(s\times (r+1),\mathbb{R}).$$

\vs{0.5cm}
A full preparation theorem for restricted log-exp-analytic power functions in $(u,z)$ is the following.

\vs{0.5cm}
{\bf2.16 Proposition}

\vs{0.1cm}
{\it
Let $e \in \mathbb{N}_0$. Let $f:X \to \mathbb{R}, (w,u,z) \mapsto f(w,u,z),$ be a restricted log-exp-analytic power function in $(u,z)$ of order at most $e$. Then there is  $r \in \mathbb{N}_0$ and a definable cell decomposition $\mathcal{C}$ of $X_{\neq 0}$ such that for every $C \in \mathcal{C}$ the function $f|_C$ is $(m+1,X)$-power-restricted $(e,r)$-prepared in $z$.}

\vs{0.4cm}
{\bf Proof}

\vs{0.1cm}
By Proposition 2.14 there is a decomposition $\mathcal{D}$ of $X$ into finitely many definable cells such that for every $D \in \mathcal{D}$ we have that $f|_D$ is $(m+1,X)$-power-restricted $e$-prepared. Fix $D \in \mathcal{D}$ and a corresponding finite set $L$ of log-analytic functions on $D$ from Definition 2.13. Let $L:=\{l_1,\ldots,l_\kappa\}$ for $\kappa \in \mathbb{N}$. By Fact 1.8 there is a decomposition $\mathcal{C}$ of $D_{\neq 0}$ into finitely many definable cells such that for every $C \in \mathcal{C}$ the functions $l_1,\ldots,l_\kappa$ are $r$-log-analytically prepared in $z$ with $C$-nice coefficent, $C$-nice base functions and common $C$-nice center. Fix $C \in \mathcal{C}$ and the corresponding center $\Theta:=(\Theta_0,\ldots,\Theta_r)$ for this preparation. We proceed by induction on $e$. For $e=-1$ there is nothing to show. 

\vs{0.2cm}
$e-1 \to e$:
We have
$$f|_C=\sigma_{b_0}b_0\exp(d_0)v(\sigma_{b_1}b_1\exp(d_1),\ldots,\sigma_{b_s}b_s\exp(d_s))$$
where $\sigma_{b_0},\ldots,\sigma_{b_s} \in \{-1,0,1\}$, $b_j=\prod_{i=1}^k h_{ji}^{\lambda_{ji}}$ with $k \in \mathbb{N}$, $\lambda_{ji} \in \mathbb{R}$ and $h_{ji}$ is a positive log-analytic function on $C$ with $h_{ji} \in L$ for $i \in \{1,\ldots,k\}$ and $j \in \{0,\ldots,s\}$. Additionally $d_0,\ldots,d_s:C \to \mathbb{R}$ are locally bounded in $z$ with reference set $X$ and are $(m,X)$-power-restricted $(e-1)$-prepared. We have that $\sigma_{b_j}b_j(t,z)\exp(d_j(t,z)) \in [-1,1]$ for $(t,z) \in C$ and the function $v:[-1,1]^s \to \mathbb{R}$ is a power series which converges absolutely on an open neighbourhood of $[-1,1]^s$ with $v([-1,1]^s) \subset \mathbb{R}_{>0}$. By the inductive hypothesis we have that 
$d_0,\ldots,d_s$ are $(m+1,X)$-power-restricted $(e-1,r)$-prepared in $z$ with center $\Theta$. Let $j \in \{0,..,s\}$. Since $h_{ji}$ is $r$-log-analytically prepared in $z$ with $C$-nice coefficient, base functions and center $\Theta$ for $i \in \{1,\ldots,k\}$ one sees immediately that
$$b_j(t,z)=\hat{a}(t)\vert{\mathcal{Y}(t,z)}\vert^{\otimes \alpha} \hat{v}(\hat{b}_1(t)\vert{\mathcal{Y}(t,z)}\vert^{\otimes p_1},\ldots, \hat{b}_{\hat{s}}(t)\vert{\mathcal{Y}(t,z)}\vert^{\otimes p_{\hat{s}}})$$
for $\alpha \in \mathbb{R}^{r+1}$, $p_i \in \mathbb{Q}^{r+1}$ and $\hat{a}:\pi(C) \to \mathbb{R}$ is a finite product of powers of $C$-nice functions, $\hat{b}_1,\ldots,\hat{b}_{\hat{s}}:\pi(C) \to \mathbb{R}$ are $C$-nice functions and $\hat{v}:[-1,1]^{\hat{s}} \to \mathbb{R}$ is a power series which converges absolutely on an open neighbourhood of $[-1,1]^{\hat{s}}$ with $\hat{v}([-1,1]^{\hat{s}}) \subset \mathbb{R}_{>0}$. We are done with composition of power series. \hfill$\blacksquare$ 

\vs{0.2cm}
With this preparation theorem we are able to prove differentiability results for restricted log-exp-analytic power functions similar as in the restricted log-exp-analytic case in \cite{8}: we have to consider preparations of restricted log-exp-analytic power functions on simple cells (see Definition 2.18 below or Definition 3.7 in \cite{8}). In \cite{8} it is shown that every $C$-nice function is log-analytic if $C$ is simple, i.e. the preparation simplifies. As in $\cite{8}$ for the restricted log-exp-analytic case and \cite{9} for the log-analytic case we call a such a preparation \textit{pure}.

\vs{0.5cm}
{\bf2.17 Definition}

\vs{0.1cm}
Let $e \in \mathbb{N}_0 \cup \{-1\}$ and $r \in \mathbb{N}_0$. By induction on $e \in \mathbb{N}_0 \cup \{-1\}$ we define that $f:C \to \mathbb{R}, (w,u,z) \mapsto f(w,u,z),$ is \textbf{purely $(m+1,X)$-power-restricted $(e,r)$-prepared in $z$} and associate a \textbf{purely preparing tuple} to this preparation. 

\vs{0.2cm}
$e=-1$: We call $f$ purely $(m+1,X)$-power-restricted $(-1,r)$-prepared in $z$ if $f$ is the zero function. A preparing tuple for $f$ is then $(0)$.

\vs{0.2cm}
$e-1 \to e$: We call $f$ purely $(m+1,X)$-power-restricted $(e,r)$-prepared in $z$ if for $(t,z) \in C$
$$f(t,z)=a(t)\vert{\mathcal{Y}(t,z)}\vert^{\otimes \alpha} \exp(d_0(t,z)) \cdot \rho(t,z)$$
where $a:\pi(C) \to \mathbb{R}$ is a finite product of powers of log-analytic functions, $\alpha \in \mathbb{R}^{r+1}$, $d_0:C \to \mathbb{R}$ is locally bounded in $(u,z)$ with reference set $X$ and is purely $(m+1,X)$-power-restricted $(e-1,r)$-prepared in $z$ and $\rho$ is a function on $C$ of the following form. There is $s \in \mathbb{N}$ such that $\rho=v \circ \phi$ where $\phi:=(\phi_1,\ldots,\phi_s):C \to [-1,1]^s$ with 
$$\phi_j(t,z)=b_j(t)\vert{\mathcal{Y}(t,z)}\vert^{\otimes \gamma_j}\exp(d_j(t,z))$$
for $j \in \{1,\ldots,s\}$ where $b_1,\ldots,b_s:\pi(C) \to \mathbb{R}$ are finite products of powers of log-analytic functions, $d_1,\ldots,d_s:C \to \mathbb{R}$ are locally bounded in $(u,z)$ with reference set $X$ and are purely $(m+1,X)$-power-restricted $(e-1,r)$-prepared in $z$, $\gamma_j:=(\gamma_{j0},\ldots,\gamma_{jr}) \in \mathbb{R}^{r+1}$ and $v$ is a power series on $[-1,1]^s$ which converges absolutely on an open neighbourhood of $[-1,1]^s$ and fulfills $v([-1,1]^s) \subset \mathbb{R}_{>0}$. A purely preparing tuple for $f$ is then 
$$(r,\mathcal{Y},a,\exp(d_0),\alpha,s,v,b,\exp(d),\Gamma)$$ 
with $b:=(b_1,\ldots,b_s)$, $\exp(d):=(\exp(d_1),\ldots,\exp(d_s))$ and
$$\Gamma:=\left(\begin{array}{cccc}
\gamma_{10}&\cdot&\cdot&\gamma_{1r}\\
\cdot&& &\cdot\\
\cdot&& &\cdot\\
\gamma_{s0}&\cdot&\cdot&\gamma_{sr}\\
\end{array}\right)\in \mathcal{M}\big(s\times (r+1),\mathbb{R}).$$

\vs{0.5cm}
{\bf2.18 Definition} (\cite{8} Definition 3.4)

\vs{0.1cm}
Let $C \subset \mathbb{R}^n \times \mathbb{R}_{\neq 0}$ be a definable cell. We call $C$ \textbf{simple} if for every $t \in \pi(C)$ we have $C_t=]0,d_t[$. 

\vs{0.5cm}
{\bf2.19 Proposition}

\vs{0.1cm}
{\it 
Let $f:X \to \mathbb{R}, (w,u,z) \mapsto f(w,u,z),$ be a restricted log-exp-analytic power function in $(u,z)$. Then there are $r \in \mathbb{N}_0$, $e \in \mathbb{N}_0 \cup \{-1\}$ and a definable cell decomposition $\mathcal{C}$ of $X$ such that for every simple $C \in \mathcal{C}$ the restriction $f|_C$ is purely $(m+1,X)$-power-restricted $(e,r)$-prepared in $z$ with center $(0)$.}

\vs{0.4cm}
{\bf Proof}

\vs{0.1cm}
By Proposition 2.16 there are $r \in \mathbb{N}_0$, $e \in \mathbb{N}_0 \cup \{-1\}$ and a definable cell decomposition $\mathcal{C}$ of $X$ such that for every $C \in \mathcal{C}$ the function $f|_C$ is $(m+1,X)$-power-restricted $(e,r)$-prepared in $z$. Fix a simple $C \in \mathcal{C}$. We show by induction on $l \in \{-1,\ldots,e\}$ that $f$ is purely $(m+1,X)$-power-restricted $(e,r)$-prepared in $z$ with center $(0)$. For $l=-1$ there is nothing to show. 

\vs{0.2cm}
$l-1 \to l$: Let 
$$(r,\mathcal{Y},a,\exp(c),\alpha,s,v,b,\exp(d),\Gamma)$$
be a preparing tuple for $f$ with $b:=(b_1,\ldots,b_s)$ and $\exp(d):=(\exp(d_1),\ldots,\exp(d_s))$. By Proposition 2.15 in \cite{5} we have that $\hat{\Theta}=0$ for every center $\hat{\Theta}$ of a $k$-logarithmic scale on $C$ (where $k \in \mathbb{N}_0$). Consequently every $C$-nice function on $\pi(C)$ is log-analytic and the center $\Theta$ of $\mathcal{Y}$ vanishes. So we have that $a$ and $b_1,\ldots,b_s$ are finite products of powers of log-analytic functions on $C$. So we see that $f$ is purely $(m+1,X)$-power-restricted $(e,r)$-prepared in $z$ with center $(0)$ by the inductive hypothesis and we are done.
\hfill$\blacksquare$

\vs{0.4cm}
A consequence of this preparation theorem is Theorem A, a version of Proposition 3.16 from \cite{8} for restricted log-exp-analytic power functions: a restricted log-exp-analytic power function $f:X \to \mathbb{R}, (w,u,z) \mapsto f(w,u,z),$ in $(u,z)$ can be real log-analytically prepared in $z$ on simple cells with coefficient and base functions which are also restricted log-exp-analytic power functions in $u$. This result is crucial for proving differentiability properties for restricted log-exp-analytic power functions.  

\vs{0.5cm}
{\bf2.20 Proposition}

\vs{0.1cm}
{\it 
Suppose that $0$ is interior point of $X_t$ for every $t \in \pi(X)$. Let $f:X \to \mathbb{R}$ be a restricted log-exp-analytic power function in $(u,z)$. Then there is $r \in \mathbb{N}_0$ and a definable cell decomposition $\mathcal{C}$ of $X$ such that for every simple $C \in \mathcal{C}$ the following holds. The restriction $f|_C$ is $r$-real log-analytically prepared with LA-preparing tuple $$(r,\mathcal{Y},a,\alpha,s,v,b,\Gamma)$$ 
where $a$ and $b_1,\ldots,b_s$ are restricted log-exp-analytic power functions in $u$ with reference set $\pi(X)$ and $\mathcal{Y}$ is an $r$-logarithmic scale with center $0$ on $C$.}

\vs{0.4cm}
{\bf Proof}

\vs{0.1cm}
The proof of this theorem is very similar to the proof of Proposition 3.16 in \cite{8}. For the readers convenience we give some details. 

\vs{0.1cm}
By Proposition 2.19 there are $r \in \mathbb{N}_0$, $e \in \mathbb{N}_0 \cup \{-1\}$ and a definable cell decomposition $\mathcal{Q}$ of $X$ such that for every simple $Q \in \mathcal{Q}$ the restriction $f|_Q$ is purely $(m+1,X)$-power-restricted $(e,r)$-prepared in $z$. Fix such a simple $Q \in \mathcal{Q}$. 
The following claim is the analogue of the corresponding claim from the proof of Proposition 3.16 from $\cite{8}$ to our situation. We omit its proof since it is the same as in $\cite{8}$ by replacing ''log-analytically prepared'' with ''real log-analytically prepared'' and ''restricted-log-exp-analytic'' with ''restricted log-exp-analytic power function''.

\vs{0.3cm}
{\bf Claim}

\vs{0.1cm}
Let $h$ be locally bounded in $(u,z)$ with reference set $X$ and $r$-real log-analytically prepared in $z$ with coefficient and base functions which are restricted log-exp-analytic power functions in $u$ with reference set $\pi(X)$. Then there is a definable simple set $D \subset Q$ with $\pi(D)=\pi(Q)$ such that $h=h_1 + h_2$ where
\begin{itemize}
	\item [(1)] $h_1:\pi(D) \to \mathbb{R}$ is a function such that $\exp(h_1):\pi(D) \to \mathbb{R}$ is a restricted log-exp-analytic power function in $u$ with reference set $\pi(X)$ and
	\item [(2)] $h_2:D \to \mathbb{R}$ is a bounded function such that $\exp(h_2)$ is $r$-real log-analytically prepared in $z$ with coefficient $1$ and base functions which are restricted log-exp-analytic power functions in $u$ with reference set $\pi(X)$.  
\end{itemize}

\vs{0.3cm}
We show by induction on $e$ that there is a simple definable $A \subset Q$ with $\pi(A)=\pi(Q)$ such that $f|_A$ is $r$-real log-analytically prepared in $z$ with coefficient and base functions which are restricted log-exp-analytic power functions in $u$ with reference set $\pi(X)$. For $l=-1$ it is clear by choosing $A:=Q$. 

\vs{0.2cm}
$l-1 \to l:$ Let 
$$(r,\mathcal{Y},a,e^d,\alpha,s,v,b,e^c,\Gamma)$$
be a purely preparing tuple for $f$ where $b:=(b_1,\ldots,b_s)$, $e^c:=(e^{c_1},\ldots,e^{c_s})$, and $\Gamma:=(\gamma_1,\ldots,\gamma_s)^t$. Note that $a,b_1,\ldots,b_s$ are finite products of powers of log-analytic functions and that $d,c_1,\ldots,c_s$ are purely $(m+1,X)$-power-restricted $(l-1,e)$-prepared in $z$. We have
$$f(t,z)=a(t)\vert{\mathcal{Y}(z)}\vert^{\otimes \alpha}e^{d(t,z)} v(b_1(t)\vert{\mathcal{Y}(z)}\vert^{\otimes \gamma_1}e^{c_1(t,z)},\ldots,b_s(t)\vert{\mathcal{Y}(x)}\vert^{\otimes \gamma_s}e^{c_s(t,z)})$$
for every $(t,z) \in Q$. By the inductive hypothesis and the claim we find a simple definable set $A \subset Q$ with $\pi(A)=\pi(Q)$ and functions $d_1,c_{11},\ldots,c_{1s}:\pi(A) \to \mathbb{R}$ and $d_2,c_{21},\ldots,c_{2s}:A \to \mathbb{R}$ with the following properties: 

\begin{itemize}
	\item [(1)] The functions $\exp(d_1)$ and $\exp(c_{11}),\ldots,\exp(c_{1s})$ are restricted log-exp-analytic power functions in $u$ with reference set $\pi(X)$,
	\item [(2)] the functions $\exp(d_2)$ and $\exp(c_{21}),\ldots,\exp(c_{2s})$ are $r$-real log-analytically prepared in $x$ with coefficient $1$ and base functions which are restricted log-exp-analytic power functions in $u$ with reference set $\pi(X)$,
	\item [(3)] we have $d|_A=d_1+d_2$ and $c_j|_A=c_{1j}+c_{2j}$ for $j \in \{1,\ldots,s\}$.
\end{itemize}

Since $a$ and $b_1,\ldots,b_s$ are products of powers of log-analytic functions we see that the functions 
$$\hat{a}:\pi(A) \to \mathbb{R}, (w,u) \mapsto a(w,u) \exp(d_1(w,u)),$$
and
$$\hat{b}_j:\pi(A) \to \mathbb{R}, (w,u) \mapsto b_j(w,u) \exp(c_{1j}(w,u)),$$
for $j \in \{1,\ldots,s\}$ are restricted log-exp-analytic power functions in $u$ with reference set $\pi(X)$. For $(w,u,z) \in A$ we have 
$$f(w,u,z)=\hat{a}(w,u)\vert{\mathcal{Y}(z)}\vert^{\otimes \alpha}e^{d_2(w,u,z)} v(\hat{\phi}_1(w,u,z),\ldots,\hat{\phi}_s(w,u,z))$$
where $\hat{\phi}_j(w,u,z):=\hat{b}_j(w,u)\vert{\mathcal{Y}(z)}\vert^{\otimes \gamma_j} e^{c_{2j}(w,u,z)}$ for $(w,u,z) \in A$ and $j \in \{1,\ldots,s\}$. By composition of power series we obtain the desired $r$-real log-analytical preparation for $h$ in $z$.  

\vs{0.2cm}
So we find a simple definable set $\hat{C} \subset Q$ with $\pi(\hat{C})=\pi(Q)$ such that $f|_{\hat{C}}$ is $r$-real log-analytically prepared in $x$ with coefficient and base functions which are restricted log-exp-analytic in $u$ with reference set $\pi(X)$. With the cell decomposition theorem applied to every such $\hat{C}$ we are done (compare with \cite{2}, Chapter 3).\hfill$\blacksquare$

\section{Differentiability Properties of Restricted  \\ 
	Log-Exp-Analytic Power Functions}

Outgoing from the preparation theorem in Proposition 2.20 we give some differentiability properties of restricted log-exp-analytic power functions. We will not give all the details in the proofs, since everything from Proposition 3.17 in \cite{8} can be formulated and proven for this class of functions in a very similar way. However we will point out the relevance of our new preparation result and explain the main differences to \cite{8}. For an example we prove the first proposition in this section completely.

\vs{0.3cm}
For this section we fix $l,m \in \mathbb{N}_0$. Let $w$ range over $\mathbb{R}^l$, $u$ over $\mathbb{R}^m$ and $z$ over $\mathbb{R}$. Let $n:=l+m$, $t:=(w,u)$, and let $\pi:\mathbb{R}^n \times \mathbb{R} \to \mathbb{R}^n, (t,z) \mapsto t,$ be the projection on the first $n$ coordinates. 

\vs{0.5cm}
{\bf3.1 Proposition}

\vs{0.1cm}
{\it
	Let $X \subset \mathbb{R}^l \times \mathbb{R}^m \times \mathbb{R}$ be definable such that $X_w$ is open for every $w \in \mathbb{R}^l$. Let $f:X \to \mathbb{R}, (w,u,z) \mapsto f(w,u,z)$, be a restricted log-exp-analytic power function in $(u,z)$. Assume that $\lim_{z \searrow 0} f(t,z) \in \mathbb{R}$ for every $t \in \pi(X)$. Then $$h:\pi(X) \to \mathbb{R}, (w,u) \mapsto \lim_{z \searrow 0}f(w,u,z),$$
	is a restricted log-exp-analytic power function in $u$.}

\vs{0.2cm}
{\bf Proof}

\vs{0.1cm}
By Proposition 2.20 there is $r \in \mathbb{N}_0$ and a definable cell decomposition $\mathcal{C}$ of $X$ such that for every simple $C \in \mathcal{C}$ the restriction $f|_C$ is $r$-real log-analytically prepared in $z$ with coefficient and base functions which are restricted log-exp-analytic power functions in $u$ with reference set $\pi(X)$. Let $C \in \mathcal{C}$ be such a simple cell. Set $g:=f|_C$ and let 
$$(r,\mathcal{Y},a,\alpha,s,v,b,\Gamma)$$
be a corresponding LA-preparing tuple for $g$. Note that $\mathcal{Y}$ has center $0$. Then 
$$g(w,u,z)=a(w,u)\vert{\mathcal{Y}(z)}\vert^{\otimes \alpha}v(b_1(w,u)\vert{\mathcal{Y}(z)}\vert^{\otimes \gamma_1},\ldots,b_s(w,u)\vert{\mathcal{Y}(z)}\vert^{\otimes \gamma_s})$$ 
for $(w,u,z) \in C$. For $j \in \{1, \ldots , s\}$ we have $\lim_{z \searrow 0} \vert{\mathcal{Y}(z)}\vert^{\otimes \gamma_j} \in \mathbb{R}$ since $b_j \neq 0$ and $b_j(u,w)\vert{\mathcal{Y}(z)}\vert^{\otimes \gamma_j} \in [-1,1]$ for $(u,w,z) \in C$. Additionally we have $\lim_{z \searrow 0} \vert{\mathcal{Y}(z)}\vert^{\otimes \alpha} \in \mathbb{R}$ if $a \neq 0$ since $u([-1,1]^s) \subset [1/c,c]$ for a constant $c>1$ and $\lim_{z \searrow 0} f(u,w,z) \in \mathbb{R}$ for $(u,w) \in \pi(X)$. Therefore we see that
$$A:\pi(C) \to \mathbb{R}, (w,u) \mapsto \lim_{z \searrow 0} a(w,u)\vert{\mathcal{Y}(z)}\vert^{\otimes \alpha},$$
and, for $j \in \{1,\ldots,s\}$, that
$$B_j:\pi(C) \to [-1,1], (w,u) \mapsto \lim_{z \searrow 0} b_j(w,u)\vert{\mathcal{Y}(z)}\vert^{\otimes \gamma_j},$$
are well-defined restricted log-exp-analytic power functions in $u$ with reference set $\pi(X)$. We obtain for $(w,u) \in \pi(C)$
$$h(w,u)=A(w,u)v(B_1(w,u),\ldots,B_s(w,u)).$$
Hence $h|_{\pi(C)}$ is a restricted log-exp-analytic power function in $u$ with reference set $\pi(X)$ by Remark 2.8. By Remark 2.9 we obtain that $h$ is a restricted log-exp-analytic power function in $u$ with reference set $\pi(X)$. \hfill$\blacksquare$

\vs{0.2cm}
Now we obtain part (1) of Theorem B.

\vs{0.3cm}
{\bf3.2 Proposition} (Closedness under taking derivatives)

\vs{0.1cm}

{\it Let $X \subset \mathbb{R}^l \times \mathbb{R}^m$ be definable such that $X_w$ is open for every $w \in \mathbb{R}^l$. Let $f:X \to \mathbb{R}, (w,u) \mapsto f(w,u),$ be a restricted log-exp-analytic power function in $u$. Let $i \in \{1,\ldots,m\}$ be such that $f$ is differentiable with respect to $u_i$ on $X$. Then $\partial f/\partial u_i$ is a restricted log-exp-analytic power function in $u$.}

\vspace{0.2cm}
{\bf Proof}

\vs{0.1cm}
The proof follows exactly the proof of Theorem A in \cite{8} with the difference that we use Remark 2.11,  Proposition 2.20 and Proposition 3.1 instead of the corresponding results from \cite{8}. \hfill$\blacksquare$

\vs{0.3cm}
The fact that a restricted log-exp-analytic power function $g:Y \to \mathbb{R}, (t,z) \mapsto g(t,z),$ in $z$, where $Y \subset \mathbb{R}^n \times\mathbb{R}$ is definable and $Y_t$ is open for every $t \in \mathbb{R}^n$, can be real log-analytically prepared in $z$ on simple definable cells gives a univariate result concerning strong quasianalyticity of $g$ in $z$ at $0$: there is $N \in \mathbb{N}$ such that if $g(t,-)$ is $C^N$ at $0$ and all derivatives of $g$ of order at most $N$ with respect to $z$ vanish at $0$ then $g(t,-)$ vanishes identically at a small interval around zero (compare with the proof of Proposition 3.19 in \cite{8} with real exponents in the log-analytical preparation instead of rational ones). With this consideration we get part (2) of Theorem B.

\vs{0.4cm}
{\bf3.3 Proposition} (Strong quasianalyticity)

\vs{0.1cm}
{\it Let $X \subset \mathbb{R}^l \times \mathbb{R}^m$ be definable such that $X_w$ is open and connected for every $w \in \mathbb{R}^n$. Let $f:X \to \mathbb{R}, (w,u) \mapsto f(w,u),$ be a restricted log-exp-analytic power function in $u$. Then there is $N \in \mathbb{N}$ with the following property. If for $w \in \mathbb{R}^l$ the function $f_w$ is $C^N$ and if there is $a \in X_w$ such that all derivatives up to order $N$ vanish in $a$ then $f_w$ vanishes identically.} 

\vspace{0.2cm}
{\bf Proof}

\vs{0.1cm}
The proof follows exactly the proof of Theorem B in \cite{8} with the difference that we use Remark 2.11 and strong quasianalyticity of restricted log-exp-analytic power functions $g:Y \to \mathbb{R}, (t,z) \mapsto g(t,z),$ in $z$ at $0$ instead of the corresponding results from \cite{8}. \hfill$\blacksquare$

\vs{0.2cm}
Another consequence of Proposition 2.20 which can be immediately shown as in \cite{8} is a univariate result concerning real analyticity of a restricted log-exp-analytic power function $f:Y \to \mathbb{R}, (t,z) \mapsto f(t,z),$ in $z$ at $0$: there is $M \in \mathbb{N}$ such that if $f(t,-)$ is $C^M$ at $0$ then $f(t,-)$ is real analytic at $0$ (compare with the proof of Proposition 3.21 in \cite{8} with real exponents in the preparation instead of rational ones). With this consideration we get part (3) of Theorem B.

\vs{0.5cm}
{\bf3.4 Theorem} (Tamm's theorem)

\vs{0.1cm}
{\it Let $X \subset \mathbb{R}^l \times \mathbb{R}^m$ be definable such that $X_w$ is open for every $w \in \mathbb{R}^l$. Let $f:X \to \mathbb{R}, (w,u) \mapsto f(w,u),$ be a restricted log-exp-analytic power function in $u$. Then there is $M \in \mathbb{N}$ such that for all $w \in \mathbb{R}^l$ if $f(w,-)$ is $C^M$ at $u$ then $f(w,-)$ is real analytic at $u$.} 

\vs{0.3cm}
{\bf Proof}

\vs{0.1cm}
The proof is the same as the proof of Proposition 3.25 in \cite{8} with the minor difference that we use Remark 2.11 and the property about real analyticity of restricted log-exp-analytic power functions at $0$ instead of the corresponding results from \cite{8}. $\hfill\blacksquare$

\vs{0.5cm}
{\bf3.5 Corollary}

\vs{0.1cm}
{\it Let $X \subset \mathbb{R}^l \times \mathbb{R}^m$ be definable such that $X_w$ is open for every $w \in \mathbb{R}^l$ and let $f:X \to \mathbb{R}, (w,u) \mapsto f(w,u),$ be a restricted log-exp-analytic power function in $u$. Then the set of all $(w,u) \in X$ such that $f(w,-)$ is real analytic at $u$ is definable.}

\vs{0.5cm}
{\bf3.6 Remark}

\vs{0.1cm}
The function
$$f:\mathbb{R} \to \mathbb{R}, u \mapsto \left\{\begin{array}{cc}
e^{-\frac{1}{u}},&u>0,\\
0,& u \leq 0,
\end{array}\right.$$
is not a restricted log-exp-analytic power function in $u$.

\vs{0.3cm}
{\bf Proof}

\vs{0.1cm}
Note that $f$ is flat at $0$, but not the zero function. So $f$ is not strong quasianalytic. Furthermore $f$ is $C^\infty$ at $0$, but not real analytic. So we see with Theorem B that $f$ is not a restricted log-exp-analytic power function in $u$. \hfill$\blacksquare$

\vs{0.4cm}

\vs{1cm}
Andre Opris\\
University of Passau\\
Faculty of Computer Science and Mathematics\\
andre.opris@uni-passau.de\\
D-94030 Germany
\end{document}